\documentclass{article}
\usepackage[american]{babel}
\usepackage{amsfonts,amsmath,amssymb}
\usepackage{xypic}
\xyoption{all}
\input{xypic}
\newdir{ >}{{}*!/-8pt/\dir{>}}
\def\R{\mathbb{R}}
\def\N{\mathbb{N}}
\def\C{\mathbb{C}}

\def\Z{\mathbb{Z}}

\def\K{\mathbb{K}} 
\def\pr{$\bf{Proof.}$\quad}
\def\fin{\hfill$\square$\\}

\def\ga{\gamma}

\def\ra{\rangle}

\newtheorem{theo}{Theorem}
\newtheorem{defi}{Definition}
\newtheorem{rem}{Remark}
\newtheorem{prop}{Proposition}
\newtheorem{cor}{Corollary}
\newtheorem{lem}{Lemma}
\begin{document}
\title{On Lie algebra crossed modules} 

\author{Friedrich Wagemann \\
        Laboratoire de Math\'ematiques Jean Leray\\
	Facult\'e des Sciences et Techniques\\
        Universit\'e de Nantes\\
	2, rue de la Houssini\`ere\\
        44322 Nantes cedex 3\\
        France\\
        tel.: ++33.2.51.12.59.57\\
        e-mail: wagemann@math.univ-nantes.fr}

\maketitle

key words: crossed modules of Lie algebras, complex simple Lie algebra, vector field Lie algebras, Godbillon-Vey cocycle 

\begin{abstract}
The goal of this article is to construct a crossed module representing the cocycle $\langle[,],\rangle$ generating $H^3({\mathfrak g};\C)$ for a simple complex Lie algebra ${\mathfrak g}$.
\end{abstract}

\section*{Introduction}
The goal of this article is to construct crossed modules for some famous third cohomology classes. 

Let ${\mathfrak g}$ be a Lie algebra and $V$ be a ${\mathfrak g}$-module. A {\it crossed module} associated to ${\mathfrak g}$ and $V$ is a homomorphism $\mu$ of Lie algebras $\mu:{\mathfrak m}\to{\mathfrak n}$ together with a compatible action of ${\mathfrak n}$ on ${\mathfrak m}$ such that $ker\,\mu=V$ and $coker\,\mu={\mathfrak g}$. According to S. Mac Lane \cite{ML}, it is M. Gerstenhaber \cite{Ger64} \cite{Ger66} to whom we can attribute the theorem stating the isomorphism between the group of equivalence classes of crossed modules associated to ${\mathfrak g}$ and $V$ and $H^3({\mathfrak g};V)$. The proof given in \cite{Ger64} is rather abstract: it relies on the characterization of $Ext$ by axioms and shows that the set of equivalence classes of exact sequences which are Yoneda products of an abelian extension and an $n$-term exact sequence satisfies the same axioms as $Ext^n$, without mentioning crossed modules. On the other hand, in \cite{Ger66} Gerstenhaber uses this isomorphism to show the above mentioned theorem about crossed modules.

Perhaps it was the level of abstraction involved which made mathematicians like better Hochschild's approach \cite{Ho1} \cite{Ho2} (due initially to S. Eilenberg and S. Mac Lane who studied the case of groups). In this approach, equivalence classes of {\it Lie algebra kernels}, i.e. homomorphisms $\psi:{\mathfrak g}\to out({\mathfrak a})$ with $Z({\mathfrak a})=V$, are shown to be isomorphic to $H^3({\mathfrak g};V)$, the class in $H^3({\mathfrak g};V)$ being an obstruction for $\psi$ to come from a general Lie algebra extension. In contrast to this obstruction point of view, our point of view is constructive; we show how crossed modules arise constructively from other algebraic objects associated to ${\mathfrak g}$ and $V$. The explicit relation between the constructive and the obstructive point of view is obscure to us in general.

There is still another approach to $3$-cohomology; indeed, in the relative setting, one fixes not only the quotient Lie algebra ${\mathfrak g}$ but all the quotient map $\pi:{\mathfrak n}\to{\mathfrak g}$. Considering relative crossed modules with fixed $V$ and $\pi$ means fixing three terms in the $4$-term exact sequence associated to a crossed module and reduces the problem to an extension type problem which is solved by Kassel-Loday in \cite{KasLod}. In order to link this approach to the general setting, one presents ${\mathfrak g}$ as a quotient of a  free Lie algebra. We think it is therefore not suited for explicitly constructing meaningful crossed modules.  

We start citing the isomorphism between the group of equivalence classes of crossed modules with kernel $V$ and cokernel ${\mathfrak g}$ and $H^3({\mathfrak g};V)$. An explicit proof can serve to have a bijection even in a topological framework, where sequences are topologically split and cocycles continuous (theorem $2$). It turns out that a class (the principal construction, shown in section $2$) of very simple crossed modules arising from an abelian extension and a short exact sequence of coefficients may be used to show surjectivity easily.

It is this same construction which we use later in the paper to unveil a crossed module representing the Godbillon-Vey cocycle \cite{Wag} (section $3.1$), one representing the cocycle $\langle[,],\rangle$ generating $H^3({\mathfrak g};\C)$ for a simple complex Lie algebra ${\mathfrak g}$ (section $4$) and the Godbillon-Vey group cocycle of the group of orientation preserving $\Z$-equivariant diffeomorphisms on the real line (section $5$).

The main theorem (theorem $4$) is the construction of a crossed module representing the cocycle $\langle[,],\rangle$. We show that this is possible using (infinite-dimensional) ${\mathfrak g}$-modules from category ${\cal O}$ (while it is well-known that this is not possible using finite dimensional ${\mathfrak g}$-modules, cf \cite{Ho2}). The proof uses spectral sequence arguments. 

We are convinced that our constructions can be adapted to many other classes of algebras and to groups in the spirit of \cite{Ger64}, \cite{Ger66}.\\ 

{\bf Acknoledgements:}\\
It is a pleasure to thank K.-H. Neeb for help on general extensions of Lie algebras and many useful discussions, to thank C. Kapoudjian for help concerning section $5$ and to thank especially S. Kumar for suggesting to use the Hochschild-Serre spectral sequence in the proof of theorem $4$. We also thank T. Pirashvili for indicating the reference \cite{BauPir}, and the referee for correcting an error.

Furthermore, we thank the Seminar S. Lie and Fields Institute in Toronto for the possibility to present preliminary versions of this work.

\section{Crossed modules and cohomology}

In this section, we present the standard material on crossed modules of Lie algebras, namely we associate to a crossed module involving a Lie algebra ${\mathfrak g}$ and a ${\mathfrak g}$-module $V$ a cohomology class $[\ga]\in H^3({\mathfrak g},V)$. The main theorem of this section is the well-known \cite{Ger66} fact that the group of equivalence classes of crossed modules involving ${\mathfrak g}$ and $V$ is isomorphic to  $H^3({\mathfrak g},V)$.

\subsection{Crossed modules}

Fix a field $\K$. Mainly we consider $\K=\R$ or $\K=\C$. 

\begin{defi}
A crossed module of Lie algebras is a homomorphism of
Lie algebras $\mu:\mathfrak{m}\to\mathfrak{n}$ together with an action $\eta$ of $\mathfrak{n}$ on $\mathfrak{m}$ by derivations such that

(a) $\mu(\eta(n)\cdot m)\,=\,[n,\mu(m)]$ for all $n\in\mathfrak{n}$
and all $m\in\mathfrak{m}$,

(b) $\eta(\mu(m))\cdot m'\,=\,[m,m']$ for all $m,m'\in\mathfrak{m}$.
\end{defi}
\begin{rem}
Later in the paper, we feel free to suppress the map $\eta$ and to write the action simply as $n\cdot m$ for $n\in{\mathfrak n}$ and $m\in{\mathfrak m}$.\\
To each crossed module of Lie algebras $\mu:\mathfrak{m}\to\mathfrak{n}$, one associates a four term exact sequence
\begin{displaymath}
0\to
V\stackrel{i}{\to}\mathfrak{m}\stackrel{\mu}{\to}\mathfrak{n}\stackrel{\pi}{\to}\mathfrak{g}\to
0
\end{displaymath}
where $ker(\mu)\,=:\,V$ and $\mathfrak{g}\,:=\,coker(\mu)$.
\end{rem}
\begin{rem}
In the framework of (infinite dimensional) locally convex Lie algebras and locally convex modules, we suppose furthermore that the sequence in remark $1$ is topologically split (i.e. all images and kernels are closed and topologically direct summands) in order to have no quotient pathologies, for example in defining $\mathfrak{g}=coker(\mu)$. We shall call such a crossed module topologically split (see appendix $A$).
\end{rem}
\begin{rem}
$\bullet$ By (a), $\mathfrak{g}$ is a Lie algebra, because $im\,(\mu)$ is an ideal.\\ 
$\bullet$ By (b), $V$ is a central Lie subalgebra of $\mathfrak{m}$, and in particular abelian.\\
$\bullet$ By (a), the action of $\mathfrak{n}$ on $\mathfrak{m}$ induces a structure of a $\mathfrak{g}$-module on $V$.\\ 
$\bullet$  Note that in general
$\mathfrak{m}$ and $\mathfrak{n}$ are not $\mathfrak{g}$-modules.
\end{rem}
\begin{defi}
Two crossed modules $\mu:{\mathfrak m}\to{\mathfrak n}$ (with action $\eta$) and $\mu':{\mathfrak m}'\to{\mathfrak n}'$ (with action $\eta'$) such that $ker(\mu)\,=\,ker(\mu')=:V$ and $coker(\mu)\,=\,coker(\mu')=:{\mathfrak g}$ are called elementary equivalent if there are morphisms of Lie algebras $\phi:{\mathfrak m}\to {\mathfrak m}'$ and $\psi:{\mathfrak n}\to {\mathfrak n}'$ such that they are compatible with the actions, meaning:
$$\phi(\eta(n)\cdot m)\,=\,\eta'(\psi(n))\cdot\phi(m)\,\,\,\,\,\,\,\forall n\in{\mathfrak n}\,\,\,\forall m\in{\mathfrak m},$$
and such that the following diagram is commutative:\\

\hspace{2cm}
\xymatrix{
0 \ar[r] & V \ar[d]^{id_V} \ar[r]^{i} & {\mathfrak m} \ar[d]^{\phi} \ar[r]^{\mu} & {\mathfrak n} \ar[d]^{\psi} \ar[r]^{\pi} &  {\mathfrak g} \ar[d]^{id_{\mathfrak g}} \ar[r] & 0 \\
0 \ar[r] & V  \ar[r]^{i'} & {\mathfrak m}' \ar[r]^{\mu'} & {\mathfrak n}'  \ar[r]^{\pi'} &  {\mathfrak g}  \ar[r] & 0}
\vspace{.5cm}

We call equivalence of crossed modules the equivalence relation generated by elementary equivalence. Let us denote by ${\rm crmod}(\mathfrak{g},V)$
the set of equivalence classes of Lie algebra crossed modules with respect to fixed kernel $V$ and fixed cokernel $\mathfrak{g}$. 
\end{defi}

\begin{rem}
In the framework of topological crossed modules, all maps in the equivalence diagram are supposed to be topologically split.
\end{rem}

\begin{theo}[Gerstenhaber]
There is an isomorphism of abelian groups
\begin{displaymath}
b:{\rm crmod}(\mathfrak{g},V)\,\cong\,H^3(\mathfrak{g},V).
\end{displaymath}
\end{theo}
Here we suppose crossed modules not necessarily split and cocycles not necessarily continuous, i.e. we are in the non-topological setting. In the topological setting, we have: 

\begin{theo}
Denote by ${\rm crmod}_{\rm top}(\mathfrak{g},V)$ the abelian group of topologically split crossed modules and by $H^3_{\rm top}(\mathfrak{g},V)$ the abelian group of continuous cohomology classes\footnote{We use indices here to emphasize the difference between the two settings - in the following, we will omit the index ``top''.}. Suppose that there is a topologically split exact sequence of $\mathfrak{g}$-modules 
$$0\to V\to W\to U\to 0$$
such that $H^3_{\rm top}({\mathfrak g},W)=0$.\\
Then there is an isomorphism of abelian groups
\begin{displaymath}
b:{\rm crmod}_{\rm top}(\mathfrak{g},V)\,\cong\,H^3_{\rm top}(\mathfrak{g},V).
\end{displaymath}
\end{theo}
An abstract proof of theorem $1$ is given in \cite{Ger66} relying on the characterization of the $Ext$-groups by axioms cf. \cite{Ger64}. Using instead explicitely cocycles gives theorem $2$ by elementary homological algebra. Let us just recall for the purpose of theorem $3$ how to associate to a crossed module a $3$-cocycle of $\mathfrak{g}$ with values in $V$: recall the exact sequence from remark $1$:
\begin{displaymath}
0\to
V\stackrel{i}{\to}\mathfrak{m}\stackrel{\mu}{\to}\mathfrak{n}\stackrel{\pi}{\to}\mathfrak{g}\to
0
\end{displaymath}
The first step is to take a linear section $\rho$ of $\pi$ and
to calculate the default of $\rho$ to be a Lie algebra homomorphism,
i.e.
\begin{displaymath}
\alpha(x_1,x_2):=[\rho(x_1),\rho(x_2)]-\rho([x_1,x_2]).
\end{displaymath}
Here, $x_1,x_2\in\mathfrak{g}$. $\alpha$ is bilinear and skewsymmetric in $x_1,x_2$. We have obviously $\pi(\alpha(x_1,x_2))=0$, because $\pi$ is a Lie
algebra homomorphism, so $\alpha(x_1,x_2)\in im(\mu)=ker(\pi)$. This
means that there exists $\beta(x_1,x_2)\in\mathfrak{m}$ such that
\begin{displaymath}
\mu(\beta(x_1,x_2))\,=\,\alpha(x_1,x_2).
\end{displaymath}
Choosing a linear section $\sigma$ on $im(\mu)$, one can choose $\beta$ as 
\begin{displaymath}
\beta(x_1,x_2)\,=\,\sigma(\alpha(x_1,x_2))
\end{displaymath}
showing that we can suppose $\beta$ bilinear and skewsymmetric in $x_1,x_2$.

Now, one computes $\mu(d^{\mathfrak{m}}\beta(x_1,x_2,x_3))=0$, where $d^{\mathfrak{m}}$ is just the formal expression of the Lie algebra coboundary, while the map $\eta\circ\rho$ is not an action in general.

 Choosing a linear section $\tau$ on $i(V)=ker(\mu)$, one can choose $\gamma$ to be $\tau\circ d^{\mathfrak{m}}\beta$ (in the obvious sense) gaining that $\gamma$ is also trilinear and skewsymmetric in $x_1,x_2$ and $x_3$. One computes that $\gamma$ is a $3$-cocycle of $\mathfrak{g}$ with values in $V$. 

In order to prove theorem $1$, one has to show that $\gamma$ is independent of the choices made, and that the resulting map 
\begin{equation}\label{b}
b:{\tt crmod}({\mathfrak g},V)\to H^3({\mathfrak g},V)
\end{equation} 
is bijective. This is rather straightforward (we use \cite{BauPir} and \cite{Nee1} in the proof of injectivity, and show surjectivity below by our methods).

\section{The principal construction}

In this section, we show how to construct a crossed module associated to a short exact sequence of $\mathfrak{g}$-modules and an abelian extension of $\mathfrak{g}$.

\subsection{Crossed modules and coefficient sequences}

Let us ask the following question: given a Lie algebra $\mathfrak{g}$, a short exact sequence of $\mathfrak{g}$-modules
\begin{equation}\label{*}
0\to V_1\to V_2\to V_3\to 0
\end{equation}
(regarded as a short exact sequence of abelian Lie algebras) and an abelian extension $\mathfrak{e}$ of $\mathfrak{g}$ by the abelian Lie algebra $V_3$
\begin{equation}\label{**}
0\to V_3\to \mathfrak{e}\to \mathfrak{g}\to 0,
\end{equation}
is the Yoneda product of (\ref{*}) and (\ref{**}) a crossed module ? In the framework of continuous cohomology, we suppose both sequences to be topologically split. In case (\ref{**}) is given by a $2$-cocycle $\alpha$, we use the notation $\mathfrak{e}=V_3\times_{\alpha}{\mathfrak g}$.
\begin{theo}
In the above situation, the Yoneda product of $(\ref{*})$ and $(\ref{**})$ is a crossed module, the associated $3$-cocycle of which is (cohomologuous to) the image of the $2$-cocycle defining the central extension $(\ref{**})$ under the connecting homomorphism in the long exact cohomology sequence associated to $(\ref{*})$.
\end{theo}
\pr
Splicing the sequences (\ref{*}) and (\ref{**}) together, one gets a map 
\begin{equation}\label{***}
\mu:V_2\to\mathfrak{e}.
\end{equation} 
Writing $\mathfrak{e}=V_2\oplus\mathfrak{g}$ as vector spaces, we have $\mu(v)=(v,0)$. On the other hand, the $\mathfrak{e}$-action $\eta$ on $V_2$ is induced by the action of $\mathfrak{g}$ on $V_2$: $\eta(w,x)\cdot v:= x\cdot v$, where $(w,x)\in\mathfrak{e}$, $v\in V_2$. With these structures, condition (b) for a Lie algebra crossed module is trivially true, while condition (a) is true by definition of the bracket in the abelian extension: $\mu(\eta(w,x)\cdot v)=(x\cdot v,0)=[(w,x),(v,0)]$.

Now let us discuss the second claim. The short exact sequence (\ref{*}) induces a short exact sequence of complexes
$$0\to C^*(\mathfrak{g},V_1)\stackrel{i}{\to} C^*(\mathfrak{g},V_2)\stackrel{\pi}{\to}C^*(\mathfrak{g},V_3)\to 0.$$
Take a cocycle $\alpha\in C^2(\mathfrak{g},V_3)$, then the connecting homomorphism $\partial$ is defined as follows:

\hspace{1cm}\xymatrix{
 & \beta\in\pi^{-1}(\alpha) \ar@{|-}[d] \ar@{|->}[r]^{\pi} & \alpha \ar@{|-}[d] \\
C^2(\mathfrak{g},V_1) \ar@{ >->}[r]^{i} \ar[d]^{d^{V_1}} & C^2(\mathfrak{g},V_2) \ar@{->>}[r]^{\pi} \ar@{-}[d]^{d^{V_2}} & C^2(\mathfrak{g},V_3) \ar@{-}[d]^{d^{V_3}} \\
C^3(\mathfrak{g},V_1) \ar@{ >->}[r]^{i}  & C^3(\mathfrak{g},V_2) \ar[d]\ar@{->>}[r]^{\pi}  & C^3(\mathfrak{g},V_3) \ar[d]\\
\partial\alpha:=i^{-1}(d^{V_2}\beta) \ar@{|->}[r]^{i} & d^{V_2}\beta \ar@{|->}[r]^{\pi} & 0=d^{V_3}\alpha}\vspace{.5cm}

Here we wrote elements on top resp. on bottom of the corresponding spaces, and denoted by $d^{V_1}$, $d^{V_2}$ and $d^{V_3}$ the Lie algebra coboundaries with values in $V_1$, $V_2$ and $V_3$ respectively. Summarizing, $\partial\alpha$ is constructed by choosing an element $\beta$ preimage of $\alpha$ under $\pi$, taking $d^{V_2}\beta$ and taking an preimage of $d^{V_2}\beta$ under $i$. It is obvious that this is exactly how we constructed the $3$-cocycle corresponding to a crossed module in section $1$. Observe that the map corresponding to $\pi$ here is the map $\mu$ in (\ref{***}) (modulo extending it by $0$ to the second factor), the map defining the crossed module. Having stated the coincidence of the two constructions, it remains to take for $\alpha\in C^2(\mathfrak{g},V_3)$ the cocycle defining the abelian extension (\ref{**}).\fin

\begin{rem}
The above theorem stays true in the framework of continuous cohomology and topologically split exact sequences, because a topologically split exact sequence induces a short exact sequence of complexes of continuous cochains, see appendix $A$.
\end{rem}

\subsection{Surjectivity in Gerstenhaber's theorem}

Let us show the surjectivity of the map $b$ in equation (\ref{b}): we have to show that given a cohomology class $[\gamma]\in H^3({\mathfrak g},V)$, there is a crossed module whose associated class is $[\gamma]$. $V$ is here some ${\mathfrak g}$-module. As the category of ${\mathfrak g}$-modules posesses enough injectives (this may be false when considering the category of locally convex ${\mathfrak g}$-modules), there is an injective ${\mathfrak g}$-module $I$ and a monomorphism $i:V\hookrightarrow I$. Consider now the short exact sequence of ${\mathfrak g}$-modules:
$$0\to V\to I\to Q\to 0$$
where $Q$ is the cokernel of $i$. As $I$ is injective, the long exact sequence in cohomology gives $H^2({\mathfrak g},Q)\cong H^3({\mathfrak g},V)$. Thus $[\gamma]$ corresponds under the isomorphism (which is induced by the connecting homomorphism) to a class $[\alpha]\in H^2({\mathfrak g},Q)$, and the principal construction applied to the above short exact sequence of ${\mathfrak g}$-modules and the abelian extension of ${\mathfrak g}$ by $Q$ using the cocycle $\alpha$ gives a crossed module whose class is a preimage under $b$ of $[\gamma]$.

\begin{cor}
Every crossed module is equivalent to one coming from the principal construction.
\end{cor}

\section{Consequences of the principal construction}

\subsection{Representing the Godbillon-Vey cocycle}

As an example of the principal construction, take $\mathfrak{g}=W_1$ the Lie algebra of formal vector fields in one formal variable $x$, and take the short exact sequence of $W_1$-modules 
$$0\to \R\to F_0\stackrel{d_{DR}}{\to}F_1\to 0,$$
where $F_{\lambda}$ is the $W_1$-module of formal $\lambda$-densities. By definition, $F_{\lambda}$ is the space of formal power series $a(x)(dx)^{\lambda}$ in $x$, such that a formal vector field $f(x)\frac{d}{dx}$ acts on $a(x)(dx)^{\lambda}$ by
$$f(x)\frac{d}{dx}\cdot a(x)(dx)^{\lambda}:=(fa'+\lambda af')(dx)^{\lambda}.$$
Observe that the action on $F_0$ corresponds to the action of $W_1$ on formal functions, the one on $F_n$ for $n\in\N$ to formal $n$-forms and the one on $F_{-1}$ to the adjoint action of $W_1$ on itself. Furthermore, $d_{DR}:F_0\to F_1$ is the formal de Rham differential. The above sequence is topologically split as shown in lemma $2$ below. 

Consider now the $2$-cocycle $\alpha$ with
$$\alpha\left(f(x)\frac{d}{dx},g(x)\frac{d}{dx}\right)\,=\,\left|\begin{array}{cc} f' & g' \\ f'' & g'' \end{array}\right|(x)(dx)^1,$$
and the 3-cocycle of Godbillon-Vey $\theta(0)$ where $\theta(x)$ is given by
$$\theta(x)\left(f(x)\frac{d}{dx},g(x)\frac{d}{dx},h(x)\frac{d}{dx}\right)\,=\,\left|\begin{array}{ccc} f & g & h \\ f' & g' & h' \\ f'' & g'' & h''\end{array}\right|(x).$$

\begin{lem}
With respect to the above mentioned coefficient sequence, we have
$$\partial\alpha\,=\,\theta(0)$$
\end{lem}
\pr

The key point is the fact that $d_{DR}\theta(x)=d^{\R}\alpha$. This is a way of saying that the Gelfand-Fuks cocycle 
$$\omega\left(f(x)\frac{d}{dx},g(x)\frac{d}{dx}\right)\,=\,\int_{S^1}\left|\begin{array}{cc} f' & g' \\ f'' & g'' \end{array}\right|(x)\,\,dx$$
which generates $H^2(Vect(S^1),\R)$ (where $Vect(S^1)$ is the Lie algebra of smooth vector fields on the circle $S^1$) is the integral over $S^1$ of the cocycle $\theta(x)$. Let us recall what we mean by this (we learned this formalism from Boris Shoikhet).\\
Let $M$ be a manifold of dimension $n$ admitting a global coordinate system depending differentiable on the base point, i.e. $M$ is parallelizable. Note $\phi_x:U\to M$ this coordinate system, $x\in M$ the base point and $U\subset \R^n$ open. $\phi_x$ induces the jet map $\phi_{Vect,x}:Vect(M)\to W_n$ which is a homomorphism of Lie algebras from the Lie algebra of smooth vector fields on $M$ to the Lie algebra of formal vector fields in $n$ formal variables (sending a field to its infinite Taylor jet at $x$). Therefore we get a map $\phi_{Vect,x}^*:C^*(W_n)\to C^*(Vect(M))$ (of complexes with trivial coefficients). One easily shows that the cocycles $\xi(x):=\phi_{Vect,x}^p(\xi)$ for a $\xi\in Z^p(W_n)$ are all cohomologuous, because $\phi_x$ depends differentiably on $x$. Thus we have $d_{DR}\xi(x)=d^{\R}\theta_1$ for some $\theta_1\in\Omega^1(M)\otimes C^{p-1}(Vect(M))$. Inductively, one gets $\theta_i\in\Omega^i(M)\otimes C^{p-i}(Vect(M))$ such that $d_{DR}\theta_i=d^{\R}\theta_{i+1}$. By Stokes' theorem, the integral $\int_{\sigma}\theta_l$ over a singular $l$-cycle $\sigma$ is a cocycle:
$$d^{\R}\int_{\sigma}\theta_l\,=\,\int_{\sigma}d^{\R}\theta_l\,=\,\int_{\sigma}d_{DR}\theta_{l-1}\,=\,0$$
The actual computation $d_{DR}\theta(x)=d^{\R}\alpha$ is straightforward and left to the reader.\\
Now to compute $\partial\alpha$, we first choose a primitive $A(f,g)$ for $\alpha(f,g)$, $f,g,h$ being formal vector fields on the line. Then we have to compute $d^{F_0}A(f,g,h)=-\sum_{\rm cycl.}A([f,g],h) + \sum_{\rm cycl.}fA(g,h)'$. The second term is just $\sum_{\rm cycl.}fA(g,h)'=\sum_{\rm cycl.}f\alpha(g,h)=\theta(x)(f,g,h)$, while the first term is a primitive (say an integral $\int_0^x$) of $d^{\R}\alpha$. But $d^{\R}\alpha=d_{DR}\theta(x)$, and so the first term is $\theta(0)(f,g,h)-\theta(x)(f,g,h)$. The upshot is that $\partial\alpha(f,g,h)=\theta(0)(f,g,h)$.\fin

Thus theorem $3$ shows that
$$0\to \R\to F_0\stackrel{d_{DR}\times 0}{\to}F_1\times_{\alpha}W_1\to W_1\to 0$$
is a crossed module representing the Godbillon-Vey cocycle.

Let us give a direct proof for the existence of a long exact sequence in cohomology in the present case; for the general case, see appendix A. 

\begin{lem}
$$0\to \R\stackrel{i}{\to} F_0\stackrel{d_{DR}}{\to}F_1\to 0$$
is split as a sequence of topological vector spaces and the corresponding sequence of complexes of continuous cochains is exact.
\end{lem}
\pr 
We consider $F_0$ and $F_1$ as locally convex topological vector spaces with the following topology: a sequence of formal power series $(\sum_{n=0}^{\infty} a^i_n x^n)_{i\in \N}$ converges if and only if all coefficient sequences $(a^i_n)_{i\in\N}$ converge in the usual topology of $\R$.  

As a continuous right inverse of $d_{DR}$, we can take term-by-term integration (with zero integration constant) of the formal series. Obviously, $i$ is a homomorphism (in the sense of topological vector spaces, i.e. an isomorphism onto its image which is given the topology of $F_0$). This shows the splitting claim. Now let us show that the corresponding sequence of complexes of continuous cochains is exact: exactness is clear on the left and in the middle. Surjectivity of $d_{DR}^*$ goes as follows: given a continuous cochain $c\in C^k(W_1,F_1)$, we have
$$c(x_1\otimes\ldots\otimes x_k)\,=\,\sum_{i=0}^{\infty}c_i(x_1\otimes\ldots\otimes x_k)x^i(dx)^{1}.$$
Continuity of $c$ means that $(x_1\otimes\ldots\otimes x_k)_n\to x_1\otimes\ldots\otimes x_k$ implies $c((x_1\otimes\ldots\otimes x_k)_n)\to c(x_1\otimes\ldots\otimes x_k)$ in the formal series topology. This in turn signifies that for all $i$, $c_i((x_1\otimes\ldots\otimes x_k)_n)\to c_i(x_1\otimes\ldots\otimes x_k)$ as real numbers. $d_{DR}^*$ is surjective as a map of complexes of all cochains, so there is a cochain $a\in C^k(W_1,F_0)$ (not necessarily continuous) such that $d_{DR}\circ a=c$. But continuity of 
$$(d_{DR}\circ a)(x_1,\ldots,x_k)=\sum_{i=0}^{\infty}ia_i(x_1,\ldots,x_k)x^{i-1}(dx)^0$$
implies evidently continuity of the $a_i$ and thus $a$ is a continuous cochain.\fin

\section{A crossed module representing the generator of $H^3({\mathfrak g},\C)$}

Let ${\mathfrak g}$ denote a (finite dimensional) simple complex Lie algebra. ${\mathfrak g}$ posesses an (up to normalization) unique bilinear invariant symmetric form $\langle,\rangle$, called the Killing form of ${\mathfrak g}$. It is well known that $\langle[,],\rangle\in C^3({\mathfrak g},\C)$ is a cocycle representing a generator of the one-dimensional $H^3({\mathfrak g},\C)$. In this section, we present the main result of the present article; we construct a crossed module representing the generator of $H^3({\mathfrak g},\C)$.

\subsection{The case of $sl_2(\C)$}

The case of $sl_2(\C)$ actually follows by restriction from the example given in section 3.1. Let us first recall some standard notation.

Call $W_1^{\rm pol}(\C)$ the Lie algebra of polynomial vector fields on $\R$ with values in $\C$. We regard $sl_2(\C)\subset W_1^{\rm pol}(\C)$ as the Lie subalgebra generated by the polynomial vector fields of degree $\leq 2$; denoting by $e_i=x^{i+1}\frac{d}{dx}$ where $x$ is the formal variable on $\R$, we have $sl_2(\C)={\rm span}_{\C}(e_{-1},e_0,e_1)$. Passing to generators $h:=2e_0$, $e:=e_1$ and $f:=-e_{-1}$ gives us the usual $sl_2(\C)$-basis.

Then by restricting the crossed module given in section 3.1, we get a crossed module
$$0\to \C\to F_0^{\rm pol}(\C)\stackrel{d_{DR}\times 0}{\to}F_1^{\rm pol}(\C)\times_{\alpha}sl_2(\C)\to sl_2(\C)\to 0,$$
where $F_{\lambda}^{\rm pol}(\C)$ for $\lambda=0,1$ denotes the polynomial density modules over $W_1^{\rm pol}(\C)$, coefficient functions having values in $\C$: $F_{\lambda}^{\rm pol}(\C)=\C[x](dx)^{\lambda}$. 

Let us identify $F_{\lambda}^{\rm pol}(\C)$ for $\lambda=0,1$ with standard modules of $sl_2(\C)$ representation theory, actually with modules of category ${\cal O}$.

Denote by ${\mathfrak h}$ the Cartan subalgebra of $sl_2(\C)$ and choose a Borel subalgebra ${\mathfrak b}\supset{\mathfrak h}$. Let ${\mathfrak n}$ be the nilpotent subalgebra of $sl_2(\C)$ such that ${\mathfrak n}\oplus{\mathfrak h}={\mathfrak b}$ as vector spaces. Let $M(\lambda)=Usl_2(\C)\otimes_{U{\mathfrak b}}\C_{\lambda}$ be the Verma module of highest weight  $\lambda\in{\mathfrak h}^*$, induced from the one ${\mathfrak b}$-dimensional module $\C_{\lambda}$ where ${\mathfrak h}$ acts according to $\lambda$ and ${\mathfrak n}$ trivially. It is well known that $M(\lambda)$ has a unique proper maximal submodule $N(\lambda)$ and the quotient $L(\lambda)$ is consequently irreducible. Let us denote by $M(\lambda)^*$, $N(\lambda)^*$ and $L(\lambda)^*$ the full duals (as $sl_2(\C)$-modules) of $M(\lambda)$, $N(\lambda)$ and $L(\lambda)$, and by $M(\lambda)^{\sharp}$, $N(\lambda)^{\sharp}$ and $L(\lambda)^{\sharp}$ the restricted duals. By definition, we mean by the {\it full dual} of a ${\mathfrak g}$-module $M$ the dual vector space $M^*$ with the ${\mathfrak g}$-action $(x\cdot f)(m):= - f(x\cdot m)$ for $x\in{\mathfrak g}$, $f\in M^*$ and $m\in M$. On the other hand, by the {\it restricted dual} of a graded ${\mathfrak g}$-module $M=\oplus M_i$, we mean the space $\oplus M_i^*$ with the same ${\mathfrak g}$-action as for the full dual. Note that our definition differs from definition 2.1.1, p. 41, in \cite{Kum} by the omission of the Cartan involution. All results we need later from \cite{Kum} transfer to our setting without modification.   

\begin{lem}
There is an isomorphism of short exact sequences of $sl_2(\C)$-modules:

\vspace{.5cm}
\hspace{2cm}
\xymatrix{
0 \ar[r] & \C \ar[d] \ar[r] & F_0^{\rm pol}(\C) \ar[d] \ar[r] & F_1^{\rm pol}(\C) \ar[d] \ar[r] & 0 \\
0 \ar[r] & L(0)^{\sharp}  \ar[r] & M(0)^{\sharp} \ar[r] & N(0)^{\sharp} \ar[r]  & 0}
\vspace{.5cm}
\end{lem}

\pr
We will show this identification by choosing appropriate bases.

$W_1^{\rm pol}(\C)$ acts on $F_0^{\rm pol}(\C)=\C[x](dx)^0$ by
$$f(x)\frac{d}{dx}\cdot g(x)(dx)^0\,=\,f(x)g'(x)(dx)^0$$
for $f,g$ polynomials.

For our purpose, an appropriate basis for $F_0^{\rm pol}(\C)$ is $f_i:=\frac{1}{i}x^i(dx)^0$ for $i>0$ and $f_0:=(dx)^0$. In this basis, the $sl_2(\C)$ generators $e,f,h$ defined above act as 
$$e\cdot f_i=(i+1)f_{i+1},\,\,\,\,\,f\cdot f_i=-(i-1)f_{i-1},\,\,\,\,\,h\cdot f_i=2if_i$$
for $i\geq 2$, $e\cdot f_1=2f_2$, $f\cdot f_1=-f_0$, $h\cdot f_1=2f_1$ and $f\cdot f_0=e\cdot f_0= h\cdot f_0=0$. 

On the other hand, $M(0)$ is isomorphic to $\C[f]$ with the action 
$$e\cdot f^i=-(i-1)if^{i-1},\,\,\,\,\,f\cdot f^i=f^{i+1},\,\,\,\,\,h\cdot f^i=-2if^i$$
for $i>0$ and  $h\cdot 1=e\cdot 1=0$, $f\cdot 1=f$. We take the restricted dual $M(0)^{\sharp}$ with the action
$$x\cdot g(m)\,=\,-g(x\cdot m)$$
for $x\in{\mathfrak g}$, $m\in M(0)$ and $g\in M(0)^{\sharp}$. Then the action on the basis $\phi_i:=(i-1)!(f^i)^*$ for $i\geq 2$, $\phi_1:=(f)^*$ and $\phi_0:=(1)^*$ where $(f^i)^*$ denotes the element of $M(0)^{\sharp}$ dual to $f^i$ reads as
$$e\cdot \phi_i =(i+1)\phi_{i+1},\,\,\,\,\,f\cdot \phi_i=-(i-1)\phi_{i-1},\,\,\,\,\,h\cdot\phi_i=2i\phi_i$$
for $i\geq 2$ and $e\cdot \phi_1=2\phi_2$, $f\cdot \phi_1=-\phi_0$, $h\cdot \phi_1=2\phi_1$ and $f\cdot \phi_0=e\cdot \phi_0= h\cdot \phi_0=0$.

A simple inspection shows that $f_i\mapsto \phi_i$ for $i\geq 0$ is an isomorphism of $sl_2(\C)$-modules. It clearly maps the submodule generated by $f_0$ on the one generated by $\phi_0$. By the five lemma, we must have an isomorphism of the quotients which shows the lemma.\fin

This lemma permits to express the above crossed module as 

$$0\to \C\to M(0)^{\sharp}\to N(0)^{\sharp}\times_{\alpha}sl_2(\C)\to sl_2(\C)\to 0.$$

In this form, it generalizes immediately to a general complex simple Lie algebra. This is what we are going to show in the next subsection.

\begin{lem}
The restriction of the Godbillon-Vey cocycle to $sl_2(\C)\subset W_1^{\rm pol}(\C)$ is a non-zero multiple of the cocycle $\langle[,],\rangle$.
\end{lem}

\pr
$\langle[e,f],h\rangle=8$ and $\theta(0)(e,f,h)=-4$.\fin

\begin{rem}
Instead of the procedure of restriction, we could have used a different approach to prove that there is a cocycle $\alpha$ such that the above crossed module represents the cocycle $\langle[,],\rangle$. Cohomology computations recalled in appendix $B$ show that $H^3(sl_2(\C),F_0^{\rm pol}(\C))=0$. This means that the map 
$$H^2(sl_2(\C),F_1^{\rm pol}(\C))\to H^3(sl_2(\C),\C)$$
is surjective which shows the claim.
\end{rem}

\subsection{The general case}

Let us first introduce the cochain $\alpha$ used for the construction of the crossed module for a general simple Lie algebra. Let $U{\mathfrak g}$ be the universal enveloping algebra of ${\mathfrak g}$, $U{\mathfrak g}^+$ the augmentation ideal and $(U{\mathfrak g}^+)^{\sharp}$ its restricted dual. Let $\{x_i\}$ be a basis of ${\mathfrak g}$ orthonormal with respect to the Killing form $\langle,\rangle$. Let us denote by $\alpha$ the 2-cochain $\alpha\in C^2({\mathfrak g}, (U{\mathfrak g})^{\sharp})$
\begin{equation} \label{****}
\alpha=\sum_{i,j,k}C_{ij}^k x_i^*\wedge x_j^*\otimes x_k^*.
\end{equation}
Here $C_{ij}^k$ denote the structure constants of ${\mathfrak g}$ in the basis $\{x_i\}$. Note that 

\begin{equation} \label{*****}
\langle[,],\rangle \,=\,\sum_{i,j,k}C_{ij}^k x_i^*\wedge x_j^*\wedge x_k^*
\end{equation}
which is easily shown by applying both sides to three basis elements. Notice further that the structure constants $C_{ij}^k$ are already antisymmetric in the three indices. 

Consider now the short exact sequence of ${\mathfrak g}$-modules
$$0\to \C^{\sharp}\stackrel{\epsilon^{\sharp}}{\to}(U{\mathfrak g})^{\sharp}\stackrel{i^{\sharp}}{\to} (U{\mathfrak g}^+)^{\sharp}\to 0$$
and its corresponding short exact sequence of complexes
$$0\to C^*({\mathfrak g},\C^{\sharp})\stackrel{\epsilon^*}{\to}C^*({\mathfrak g},(U{\mathfrak g})^{\sharp})\stackrel{i^*}{\to}  C^*({\mathfrak g},(U{\mathfrak g}^+)^{\sharp})\to 0$$
Our next goal is to show that $i^*(\alpha)$ for $\alpha$ given in equation (\ref{****}) is a cocycle, and then to compute $\widetilde{\partial}(i^*(\alpha))$ where $\widetilde{\partial}$ is the connecting homomorphism corresponding to the short exact sequence of complexes.

\begin{lem}
$$d^{U{\mathfrak g}^{\sharp}}\alpha=\epsilon^*(\langle[,],\rangle)$$
In particular, $i^*(\alpha)$ is a 2-cocycle. 
\end{lem}

\pr
$$d^{U{\mathfrak g}^{\sharp}}\alpha(x_1,x_2,x_3)=\sum_{\rm cycl.}\alpha([x_i,x_j],x_k) - \sum_{\rm cycl.}x_i\cdot\alpha(x_j,x_k)$$
Here the first sum is zero. Indeed, we can evalue it on an $u\in{\mathfrak g}\subset U{\mathfrak g}$ and get by equation (\ref{*****})
$$\sum_{\rm cycl.}\langle[[x_i,x_j],x_k],u\rangle$$
which is zero by the Jacobi identity.

$x_l\cdot\alpha(x_m,x_n)$ evaluated on $u\in U{\mathfrak g}$ means explicitely 
$$\sum_{i,j,k}C_{ij}^k x_i^*(x_m)\wedge x_j^*(x_n)\wedge x_k^*(x_l\cdot u)$$ 
Thus it forces on the one hand $d^{U{\mathfrak g}^{\sharp}}\alpha$ to be non-zero at most on $\C\subset U{\mathfrak g}$. On the other hand, the evaluation on $u=1$ gives $C_{mn}^l$, which means that the second sum is equal to $\epsilon^*(\langle[,],\rangle).$\fin 

\begin{cor}
$$\widetilde{\partial}(i^*(\alpha))=\langle[,],\rangle$$
\end{cor}

Choose a Cartan subalgebra ${\mathfrak h}$ of ${\mathfrak g}$. Let $M(0)^{\sharp}$, $N(0)^{\sharp}$ and $L(0)^{\sharp}$ denote the restricted duals of the Verma module $M(0)$ of highest weight $0\in{\mathfrak h}^*$, resp. of its maximal proper submodule $N(0)$, resp. of its simple quotient $L(0)\cong \C$ as in the previous subsection. The Verma module $M(\lambda)$ can also be constructed as a quotient of $U{\mathfrak g}$ by an ideal $I(\lambda)$. Consider the diagramm of ${\mathfrak g}$-modules:

\vspace{.5cm}
\hspace{2cm}
\xymatrix{
 & 0 \ar[d] & 0 \ar[d] &  \\
0 \ar[r] & L(0)^{\sharp} \ar[r] \ar[d] & \C \ar[r] \ar[d] & 0 \ar[d] \\ 
0 \ar[r] & M(0)^{\sharp} \ar[r] \ar[d] & U{\mathfrak g}^{\sharp} \ar[r] \ar[d] & I(0)^{\sharp} \ar[d] \ar[r] & 0 \\
0 \ar[r] & N(0)^{\sharp} \ar[r] \ar[d] & (U{\mathfrak g}^+)^{\sharp} \ar[r] \ar[d] & I(0)^{\sharp} \ar[d] \ar[r] & 0 \\ 
& 0  & 0  & 0  
}\vspace{.5cm}

It induces a diagramm in cohomology a part of which looks like:

\vspace{.5cm}
\hspace{2cm}
\xymatrix{
 H^2({\mathfrak g},M(0)^{\sharp}) \ar[r] \ar[d] & H^2({\mathfrak g},(U{\mathfrak g})^{\sharp}) \ar[r] \ar[d] & H^2({\mathfrak g},I(0)^{\sharp}) \ar[d]^{\cong}  \\
 H^2({\mathfrak g},N(0)^{\sharp}) \ar[r]^{\zeta} \ar[d]^{\partial} & H^2({\mathfrak g},(U{\mathfrak g}^+)^{\sharp}) \ar[r] \ar[d]^{\widetilde{\partial}} & H^2({\mathfrak g},I(0)^{\sharp}) \ar[d]  \\
 H^3({\mathfrak g},L(0)^{\sharp}) \ar[r]^{\cong} \ar[d]^0 &  H^3({\mathfrak g},\C) \ar[r] \ar[d] & 0 \ar[d] \\ 
 H^3({\mathfrak g},M(0)^{\sharp}) \ar[r]  &  H^3({\mathfrak g},U{\mathfrak g}^{\sharp}) \ar[r]  &  H^3({\mathfrak g},I(0)^{\sharp}) 
}\vspace{.5cm}

\noindent $i^*(\alpha)$ (cf equation (\ref{****}) and lemma 9) is in $H^2({\mathfrak g},(U{\mathfrak g}^+)^{\sharp})$ and $[\widetilde{\partial}i^*(\alpha)]\in H^3({\mathfrak g},\C)$ is the wanted class. By abuse of notation, we will write $\alpha$ for $i^*(\alpha)$ from now on. Thus we have on the one hand:

\begin{prop} 
A crossed module representing $\langle[,],\rangle$ for a general simple complex Lie algebra ${\mathfrak g}$ is:\\
$$0\to \C\to U{\mathfrak g}^{\sharp}\to (U{\mathfrak g}^+)^{\sharp}\times_{\alpha}{\mathfrak g}\to {\mathfrak g}\to 0$$
where $\alpha$ is given by
$$\alpha\,=\,\sum_{i,j,k}C_{ij}^k\,\,x_i^*\wedge x_j^*\otimes x_k^*$$
\end{prop}

This proposition is not completely satisfactory as $(U{\mathfrak g}^+)^{\sharp}$ and $(U{\mathfrak g})^{\sharp}$ are not modules of category ${\cal O}$ and as we wanted the crossed module to be induced from the short exact sequence involving $M(0)^{\sharp}$, $N(0)^{\sharp}$ and $L(0)^{\sharp}$. 

On the other hand, we will show below that $H^2({\mathfrak g},L(0)^{\sharp})\to H^3({\mathfrak g},M(0)^{\sharp})$ is the zero map. This shows that $\alpha$ represents a class in $H^2({\mathfrak g},N(0)^{\sharp})$ (the map $\zeta$ being induced by an inclusion of a submodule). Certainly, $H^i({\mathfrak g},U{\mathfrak g}^*)=0$ for all $i>0$ as the full dual $U{\mathfrak g}^*$ is an injective ${\mathfrak g}$-module. $U{\mathfrak g}$ is still a $({\mathfrak g},{\mathfrak h})$-module (\cite{Kum} definition 3.3.4, p. 91), and lemma 3.3.5 {\it loc. cit.} shows that we also have $H^i({\mathfrak g},U{\mathfrak g}^{\sharp})=0$ for all $i>0$ for the restricted dual $U{\mathfrak g}^{\sharp}$. The upshot of this reasoning is that the class of $\alpha$ (under the given identifications) is also a preimage under $\partial$ of the class of $\langle[,],\rangle$. Thus we found the $\alpha$ which gives the wanted crossed module:

\begin{theo}
The short exact sequence
$$0\to L(0)^{\sharp}\to M(0)^{\sharp}\to N(0)^{\sharp}\to 0$$
and the $2$-cocycle $\alpha$ given by equation (\ref{****}) give via the principal construction a crossed module which represents (a non-zero multiple of) the class of $\langle[,],\rangle$ in $H^3({\mathfrak g},\C)$.
\end{theo}

\begin{rem}
Before we proceed to the proof of the theorem, let us remark that the above big diagramm can be resumed in the following exact sequences for $k\geq 1$, recalling the fact that $H^*({\mathfrak g},\C)$ lives only in odd degrees:

\vspace{0.2cm}
\hspace{-2.5cm}
\begin{parbox}{16cm}
{$$\scriptstyle 0\to\begin{array}{c}\scriptstyle H^{2k}({\mathfrak g},M(0)^{\sharp}) \\ \scriptstyle|| \\\scriptstyle H^{2k-1}({\mathfrak g},I(0)^{\sharp})\end{array}\to H^{2k}({\mathfrak g},N(0)^{\sharp})\to\begin{array}{c}\scriptstyle H^{2k+1}({\mathfrak g},L(0)^{\sharp}) \\\scriptstyle || \\\scriptstyle H^{2k}({\mathfrak g},(U{\mathfrak g}^+)^{\sharp})  \end{array}\to \begin{array}{c}\scriptstyle H^{2k+1}({\mathfrak g},M(0)^{\sharp}) \\\scriptstyle || \\\scriptstyle H^{2k}({\mathfrak g},I(0)^{\sharp})\end{array}\to  H^{2k+1}({\mathfrak g},N(0)^{\sharp})\to 0$$}
Here the middle term is either $\C$ or $0$. Consider also the same diagram and sequences for the case of general $\lambda\in{\mathfrak h}^*$. Here the middle term is always zero. By the cohomology results of $\cite{Wil}$, it should be possible to calculate all the terms. Observe also that the proof of theorem $4$ shows that the arrow $H^{3}({\mathfrak g},L(0)^{\sharp})\to H^3({\mathfrak g},M(0)^{\sharp})$ is zero, so the sequence splits into two for $k=1$. Does it split for $k\geq 2$ ?
\end{parbox}
\end{rem}

\noindent{\bf Proof of theorem 4.}\\
By what we have said above, the proof of the theorem reduces to the proof of the fact that $\phi:H^3({\mathfrak g},L(0)^{\sharp})\to H^3({\mathfrak g},M(0)^{\sharp})$ is the zero map. We will show this by a spectral sequence argument. Using lemma 3.3.5 of \cite{Kum}, we will no longer distinguish between cohomology with values in full duals or restricted duals. 

The strategy is to understand the map $\phi$ as a map between Hochschild-Serre spectral sequences corresponding to the subalgebra ${\mathfrak h}\subset{\mathfrak g}$. For this purpose, we denote by $E_r^{p,q}(L)$ the Hochschild-Serre spectral sequences corresponding to the filtration of $C^*({\mathfrak g},L(0)^{\sharp})$ with respect to the subalgebra ${\mathfrak h}\subset{\mathfrak g}$, and by $E_r^{p,q}(M)$ the same spectral sequence arising from filtration of $C^*({\mathfrak g},M(0)^{\sharp})$ . The map $\phi:C^*({\mathfrak g},L(0)^{\sharp})\to C^*({\mathfrak g},M(0)^{\sharp})$ maps just a cochain of some filtration with values in the submodule $L(0)^{\sharp}$ to the same cochain with the same filtration, but considered as cochain with values in the ambient module $M(0)^{\sharp}$. Thus it induces a map $\phi_r^{p,q}:E_r^{p,q}(L)\to E_r^{p,q}(M)$ for all $r\geq 0$. To show that the map $\phi$ is zero, we will show that in $\bigoplus_{p+q=3}E_r^{p,q}(L)$ only $E^{2,1}(L)$ survives to infinity, whereas in the sequence $E_r^{p,q}(M)$, $E^{2,1}(M)$ does not survive. In conclusion, $\phi:\bigoplus_{p+q=3}E_{\infty}^{p,q}(L)\to \bigoplus_{p+q=3}E_{\infty}^{p,q}(M)$ is zero and this shows the theorem.

Let us recall the Hochschild-Serre spectral sequence on which we base all of our reasoning, cf \cite{Kum} theorem E.13:

\begin{theo}
Let ${\mathfrak g}$ be a finite dimensional Lie algebra, ${\mathfrak h}\subset{\mathfrak g}$ be a reductive subalgebra and $M$ be a ${\mathfrak g}$-module which is the sum of its finite-dimensional irreducible ${\mathfrak h}$-submodules. 
Then there exists a convergent cohomology spectral sequence with
$$E^{p,q}_2\cong H^p({\mathfrak g},{\mathfrak h};M)\otimes H^q({\mathfrak h};\C)$$
and converging to $H^{p+q}({\mathfrak g};M)$.
\end{theo}

Now we discuss $E_r^{p,q}(L)$. For this we use proposition 3.2.11 in \cite{Kum}. It states that $H^i({\mathfrak g},{\mathfrak h};\C)=0$ for $i$ odd and that $H^0({\mathfrak g},{\mathfrak h};\C)\cong\C$ and $H^2({\mathfrak g},{\mathfrak h};\C)\cong{\mathfrak h}^*$ as vector spaces. From this we draw the $E_2$-plane:

\vspace{.5cm}
\hspace{2cm}
\xymatrix{
\bullet\Lambda^3{\mathfrak h}^* \ar[rrd]^{d_2^{0,3}} & \bullet 0 & \bullet\Lambda^3{\mathfrak h}^*\otimes{\mathfrak h}^* & \bullet 0  \\
\bullet\Lambda^2{\mathfrak h}^*  & \bullet 0 & \bullet\Lambda^2{\mathfrak h}^*\otimes{\mathfrak h}^* & \bullet 0 \\
\bullet{\mathfrak h}^*  & \bullet 0 & \bullet{\mathfrak h}^*\otimes{\mathfrak h}^* & \bullet 0 \\
\bullet\C  & \bullet 0 & \bullet{\mathfrak h}^* & \bullet 0 \\
}\vspace{.5cm}

\begin{lem}
The differential 
$$d_2^{0,3}:\Lambda^3{\mathfrak h}^*\cong H^0({\mathfrak g},{\mathfrak h};\C)\otimes H^3({\mathfrak h},\C)\to {\mathfrak h}^*\otimes\Lambda^2{\mathfrak h}^*\cong H^2({\mathfrak g},{\mathfrak h};\C)\otimes H^2({\mathfrak h},\C)$$
is injective.
\end{lem}

\pr Given a $c\in C^3({\mathfrak h},\C)$, we have to show that there exist $x_1,x_2\in {\mathfrak g}/{\mathfrak h}$ and $y_1,y_2\in{\mathfrak h}$ such that $d_2^{0,3}c(x_1,x_2,y_1,y_2)\not=0$. But  $d_2^{0,3}c(x_1,x_2,y_1,y_2)$ is just the usual Chevalley-Eilenberg coboundary operator applied to $c$, seen as a cochain on ${\mathfrak g}$ being non-zero only on ${\mathfrak h}$. This coboundary has six terms; the only possibly non-zero term is $c([x_1,x_2],y_1,y_2)$. In this case $deg(x_1)=-deg(x_2)$, $[x_1,x_2]\in{\mathfrak h}$. Now for a given $c$, it is clear that one can always find $x_1,x_2,y_1,y_2$ such that this term is non-zero.\fin

\begin{cor}
The only term which survives to infinity in $\bigoplus_{p+q=3}E^{p,q}_r(L)$ is $E^{2,1}_{\infty}(L)$.
\end{cor}

Now let us turn to the spectral sequence $E_r^{p,q}(L)$. By theorem $5$ and exercise 9.3.E (3) in \cite{Kum}, we can draw the $E_2$-plane:

\vspace{.5cm}
\hspace{2cm}
\xymatrix{
\bullet\Lambda^3{\mathfrak h}^*  & \bullet \Lambda^3{\mathfrak h}^* & \bullet \Lambda^3{\mathfrak h}^* & \bullet \Lambda^3{\mathfrak h}^*  \\
\bullet\Lambda^2{\mathfrak h}^* \ar[rrd]^{d_2^{0,2}} & \bullet \Lambda^2{\mathfrak h}^* & \bullet\Lambda^2{\mathfrak h}^* & \bullet \Lambda^2{\mathfrak h}^*  \\
\bullet{\mathfrak h}^* \ar[rrd]^{d_2^{0,1}}  & \bullet {\mathfrak h}^* \ar[rrd]^{d_2^{1,1}} & \bullet{\mathfrak h}^* & \bullet {\mathfrak h}^*  \\
\bullet\C  & \bullet \C & \bullet\C & \bullet \C \\
}\vspace{.5cm}

The main step towards proving that here $E_{\infty}^{2,1}(M)=0$ is the following lemma:

\begin{lem}
The differential 
$$d_2^{1,1}:{\mathfrak h}^*\cong H^1({\mathfrak g},{\mathfrak h};M(0)^{\sharp})\otimes H^1({\mathfrak h},\C)\to \C\cong H^3({\mathfrak g},{\mathfrak h};M(0)^{\sharp})\otimes H^0({\mathfrak h},\C)$$
is zero.
\end{lem}

\pr
Using the root system associated to our choice of the Cartan subalgebra ${\mathfrak h}\subset{\mathfrak g}$, ${\mathfrak g}=\bigoplus_{\alpha\in\triangle}{\mathfrak g}_{\alpha}\oplus{\mathfrak h}$ is a graded Lie algebra and $M(0)^{\sharp}$ is a graded ${\mathfrak g}$-module. These gradings induce gradings on the space of cochains:
$$C^q_{(\lambda)}({\mathfrak g},M(0)^{\sharp})=\{c\in C^q({\mathfrak g},M(0)^{\sharp})\,|\,c(x_1,\ldots,x_q)\in M(0)^{\sharp}_{(\lambda_1+...+\lambda_q-\lambda)}\,\,\,{\rm for}\,\,\,x_i\in{\mathfrak g}_{\lambda_i}\}$$
It is well known (cf \cite{Fuk} theorem 1.5.2 b) that the inclusion of the subcomplex of cochains of degree 0 induces an isomorphism in cohomology.

Let us thus restrict the cochain complex $C^*({\mathfrak g},M(0)^{\sharp})$ to the subcomplex $C^*_{(0)}({\mathfrak g},M(0)^{\sharp})$ consisting of cochains of degree 0. In this subcomplex, we construct the Hochschild-Serre spectral sequence with respect to the subalgebra ${\mathfrak h}$. As before, the second term of the spectral sequence is $E_2^{p,q}=H^p({\mathfrak g},{\mathfrak h};M(0)^{\sharp})\otimes H^q({\mathfrak h},\C)$, but here all representatives of classes viewed as cochains in $C^*({\mathfrak g},M(0)^{\sharp})$ are cochains are of degree zero. Indeed, by definition of $c\in C^p({\mathfrak g},{\mathfrak h};M(0)^{\sharp})$, we have $i_xc=0$ for all $x\in{\mathfrak h}$ and ${\cal L}_xc=0$ for all $x\in{\mathfrak h}$. This last condition means that the cochains $c$ is of degree zero.\\
Now compute 
$$d_2^{1,1}:{\mathfrak h}^*\cong H^1({\mathfrak g},{\mathfrak h};M(0)^{\sharp})\otimes H^1({\mathfrak h},\C)\to \C\cong H^3({\mathfrak g},{\mathfrak h};M(0)^{\sharp})\otimes H^0({\mathfrak h},\C).$$
$$d_2^{1,1}c(x_1,x_2,x_3)=\sum(-1)^{\sharp}c([x_i,x_j],...)+\sum(-1)^{\sharp}x_i\cdot c(...).$$
Here $c\in C^1({\mathfrak g},{\mathfrak h};M(0)^{\sharp})\otimes C^1({\mathfrak h},\C)$ and $x_i\in {\mathfrak g}/{\mathfrak h}$. The second term in this sum is zero, because there are two elements of ${\mathfrak g}/{\mathfrak h}$ as arguments of $c$. Actually, the first term is also zero: in order to have a non-zero contribution, there must be $x_1,x_2\in{\mathfrak g}/{\mathfrak h}$ with $[x_i,x_j]\in{\mathfrak h}$ on which $c$ is non-zero. This means $deg(x_1)=-deg(x_2)$. But then in order to have $deg(c)=0$, as $c$ has values in $1\otimes\C_0$ when the ${\mathfrak g}/{\mathfrak h}$-argument is in ${\mathfrak n}$, we must have $deg(x_3)=0$ which means $x_3=0$. This shows the lemma.\fin

\begin{cor}
$E_{\infty}^{2,1}(M)=0$
\end{cor}

\pr
The strategy is to see which terms survive, knowing that the final outcome must be 
$$H^p({\mathfrak g},{\mathfrak h};M(0)^{\sharp})=\Lambda^p{\mathfrak h}^*$$
for all $p\geq 0$ which follows from theorem 4.15 in \cite{Wil}.

The term $E_2^{0,0}(M)$ certainly survives to infinity and gives $H^0({\mathfrak g},{\mathfrak h};M(0)^{\sharp})=\C$. The term $E_2^{1,0}(M)$ also survives, and in order to have $H^1({\mathfrak g},{\mathfrak h};M(0)^{\sharp})={\mathfrak h}^*$, $d_2^{0,1}$ must be surjective. This means on the other hand that the term $E^{0,2}_{\infty}(M)=0$. By the above lemma, $E^{1,1}_2(M)={\mathfrak h}^*$ survives to infinity. 

Notice that the differential $d_3^{0,2}$ is zero: indeed, for a representative cochain $c\in\Lambda^3{\mathfrak h}^*$, we have 
$$d_3^{0,2}c(x_1,x_2,x_3)=\sum_{\rm cycl.}c([x_i,x_j],x_k)$$
with $x_i\in{\mathfrak g}/{\mathfrak h}$. Therefore $E^{0,2}_2(M)$ can loose dimensions only via $d_2^{0,2}$.

As a consequence, in order to have then $H^2({\mathfrak g},{\mathfrak h};M(0)^{\sharp})=\Lambda^2{\mathfrak h}^*$, $d_2^{0,2}$ must be surjective, which means exactly that $E_2^{2,1}(M)$ does not survive to infinity.\fin

This ends the proof of theorem $4$.\fin

\section{A crossed module for diffeomorphism groups}

Let $Diff(S^1)$ denote the group of orientation preserving diffeomorphisms of the circle. It is a Fr\'echet manifold, homotopy equivalent to $S^1$, and has as its universal covering group the group $\widetilde{Diff(S^1)}=Diff_{\Z}(\R)$ of orientation preserving $\Z$-equivariant diffeomorphisms of $\R$. The covering sequence reads:
$$0\to \Z\to \widetilde{Diff(S^1)}\stackrel{\pi}{\to} Diff(S^1)\to 1$$

In this section, we present a group crossed module for the group $\widetilde{Diff(S^1)}$ which correspopnds to the Godbillon-Vey cocycle. In order to explain its origin, we first discuss the construction (given in \cite{Wag}) of a crossed module representing the Godbillon-Vey cocycle for $Vect(S^1)$, the Lie algebra of smooth vector fields on $S^1$.

By the example of section $3.1$, we have a construction like this for $W_1$. One cannot transpose it directly to $Vect(S^1)$: denote by ${\cal F}_{\lambda}$ the $\lambda$-density $Vect(S^1)$-module consisting as a vector space of smooth functions on $S^1$. Now the corresponding sequence of $\lambda$-density modules reads
$$0\to\R\to{\cal F}_0\to{\cal F}_1\to H^1(S^1,\R)=\R\to 0.$$
Obviously, the last term prevents us from using the previous construction.

The trick is to pass in an appropriate manner to the universal covering: lifting elements of $Vect(S^1)$ to $2\pi$-periodic functions on $\R$, we can make them act on $F_0$ and $F_1$, the modules of densities on $\R$.

Now let us pass to group level. The relation of continuous group cohomology, cohomology as a topological space and Lie algebra cohomology of the Lie algebra of a Lie group is given by the {\it Hochschild-Mostow spectral sequence} (see \cite{Fuk} p. 295 thm. 3.4.1). One can refine this sequence to start from the cohomology of the homogeneous space of the group by some compact subgroup, abutting then to the relative Lie algebra cohomology modulo this subgroup (see \cite{Fuk} p. 297 thm. 3.4.2). An isomorphism between continuous group and relative Lie algebra cohomology stemming from acyclicity of the homogeneous space is then known as {\it van Est's theorem} (cf. \cite{Fuk} p. 298 corollary). 

The diffeomorphism group version of van Est's theorem (\cite{BroTra} thm. 1 for the version with trivial coefficients which can be adapted to coefficients cf. \cite{OvsRog} p. 286) with non-trivial coefficients implies that:

\begin{theo}
$$H^*_c(Diff(S^1);\underline{{\cal F}_{\lambda}})\,\cong\,H^*(Vect(S^1),SO(2);{\cal F}_{\lambda})$$
\end{theo}

Here $H_c^*(Diff(S^1);\underline{M})$ denotes group cohomology with continuous cochains and $\underline{M}$ is the corresponding group module for a Lie algebra module $M$. The isomorphism is functorial.

By the version with trivial coefficients (cf the corollary on p. 298 in \cite{Fuk}), we have an isomorphism
$$H^*(Diff(S^1),\R)\cong H^*(Vect(S^1),SO(2);\R)$$
and applied to $\widetilde{Diff(S^1)}$, we get an isomorphism
$$H^*(\widetilde{Diff(S^1)},\R)\cong H^*(Vect(S^1);\R).$$
The last statement follows using the fact that $\widetilde{Diff(S^1)}$ is contractible in the Hochschild-Mostow spectral sequence (see \cite{Fuk} p.295 thm. 3.4.1).

Therefore the Godbillon-Vey cocycle $\theta(0)$ arises in $H^*(\widetilde{Diff(S^1)},\R)$, but not in $H^*(Diff(S^1),\R)$. But its class is sent to the Bott-Thurston class $[\mu]$ in the Gysin sequence
$$H^3(\widetilde{Diff(S^1)},\R)\to H^2(Diff(S^1),\R),\,\,\,\,[\theta(0)]\mapsto[\mu].$$

It is easy to see that in group cohomology crossed modules can be constructed in the same way as we showed in section $2$. Moreover, theorem $3$ stays true (cf \cite{Wei} p. 188). Let us thus consider the crossed module

\begin{equation}  \label{*******}
0\to\R\to\underline{F_0}\to\underline{F_1}\times_{\pi^*B_1}\widetilde{Diff(S^1)}\to \widetilde{Diff(S^1)}\to 1.
\end{equation}

Here, $\pi^*:H^*(Diff(S^1),\underline{F_1})\to H^*(\widetilde{Diff(S^1)},\underline{F_1})$ is the map induced in cohomology by the covering projection and $B_1$ is the generator of $H^2(Diff(S^1),\underline{F_1})$ given by Ovsienko and Roger, cf theorem $9$ in appendix B. $B_1$ corresponds to the cocycle we named $\alpha$ in section $3.1$ and $\omega_1$ in theorem $8$. 

The upshot of the discussion is the following theorem:

\begin{theo}
The crossed module $(\ref{*******})$ represents the Gobillon-Vey cocycle in $H^3(\widetilde{Diff(S^1)},\R)$.
\end{theo}

\pr
The proof of the theorem follows from functoriality of the van Est isomorphism and theorem $3$ for group crossed modules. Indeed, by these two facts, we have a commutative diagramm:

\vspace{.5cm}
\hspace{2cm}
\xymatrix{
[\alpha] \ar@{|-}[d] \ar@{|->}[r] & [\pi^*B_1] \ar@{|-}[d] \\
H^2(Vect(S^1),F_1) \ar@{-}[d]^{\partial} \ar[r]^{\cong} &  H^2(\widetilde{Diff(S^1)},\underline{F_1}) \ar@{-}[d]^{\partial} \\
H^3(Vect(S^1),\R) \ar[d] \ar[r]^{\cong} & H^3(\widetilde{Diff(S^1)},\underline{\R}) \ar[d]\\
[\partial\alpha]=[\theta(0)]  \ar@{|->}[r]  & [\partial\pi^*B_1]
}\vspace{.5cm}

This shows the claim.\fin

\begin{appendix}
\section{Continuous cohomology}

In this appendix, we consider continous Lie algebra cohomology $H^*({\mathfrak g},M)$, i.e. for a topological Lie algebra ${\mathfrak g}$, we consider the Chevalley-Eilenberg complex with continuous (instead of all) cochains, with values in a topological ${\mathfrak g}$-module $M$. The vector space $H^*({\mathfrak g},M)$ is then defined as usual, cf \cite{Fuk}. The aim of this appendix is to show that a topologically split short exact sequence of modules induces a long exact sequence in continuous cohomology.  

\begin{defi}
Let $\mathfrak{g}$ be a real Fr\'echet nuclear Lie algebra, and denote by $\mathfrak{g}-{\tt mod}_{\tt top}$ the category of locally convex $\mathfrak{g}$-modules. Morphisms in this category are the continuous $\mathfrak{g}$-module homomorphisms.

A short exact sequence of objects in $\mathfrak{g}-{\tt mod}_{\tt top}$
$$0\to A\stackrel{i}{\to} B\stackrel{\pi}{\to} C\to 0$$
is called topologically split if there exists a continuous linear section $\sigma:C\to B$ of $\pi$, i.e. $\pi\circ\sigma=id_C$, and if $i$ is open onto its image. In this case, $i(A)$ has a complement in $B$: 
$$B\cong i(A)\oplus\sigma(C)$$
as topological vector spaces (in general, not as $\mathfrak{g}$-modules !).
\end{defi}

By the way, for a complex Lie algebra $\mathfrak{g}$, everything referring to $\R$ can be replaced by a reference to $\C$ without problem.

\begin{lem}
Let 
$$0\to E_1\stackrel{f}{\to} E_2\stackrel{g}{\to} E_3\to 0$$
be a topologically split short exact sequence of topological $\mathfrak{g}$-modules. Then the short exact sequence of complexes of continuous cochains
$$0\to C^p(\mathfrak{g},E_1)\stackrel{f_*}{\to} C^p(\mathfrak{g},E_2)\stackrel{g_*}{\to} C^p(\mathfrak{g},E_3)\to 0$$
is exact for all $p$. Therefore, we have a long exact sequence in continuous cohomology.
\end{lem}

\pr $f_*$ is clearly injective and $g_*\circ f_*=0$.\\
 Let us show that $ker\,g_*\subset im\,f_*$: Indeed, let $c\in ker\,g_*$, i.e. $c\in C^p(\mathfrak{g},E_2)$ such that $g\circ c=0$. Then as the coefficient sequence is topologically split, $im\,c\subset ker\,g=im\,f\cong E_1$, thus $c\in C^p(\mathfrak{g},E_1)=im\,f$.\\
Let us show that $g_*$ is surjective. Let $c\in C^p(\mathfrak{g},E_3)$. But there is a continuous linear section $\sigma$ of $g$. Then $\sigma\circ c$ is the element of $C^p(\mathfrak{g},E_2)$ we are looking for. Indeed, $g_*(\sigma\circ c)=g\circ\sigma\circ c = c$.\fin

\section{Cohomological results on $W_1$ and $sl_2(\C)$}

In this section, we recall and regroup some cohomological results (mostly well known) on $Vect(S^1)$, $W_1^{\rm pol}(\C)$ and $sl_2(\C)\subset W_1^{\rm pol}(\C)$. Notations are as in the main part of the paper.
\begin{theo}
The cohomology groups $H^p(sl_2(\C);F_{\lambda}^{\rm pol}(\C))$ are non zero only for $\lambda=0,1$, and have the following dimensions:
$$H^p(sl_2(\C),F_0^{\rm pol}(\C))\,\cong\,\left\{\begin{array}{ccc} \C^2 &{\rm for}&p=1 \\
\C &{\rm for}&p=0,2 \\ 0 &{\rm for}&p\not=0,1,2 \end{array}\right.$$
$$H^p(sl_2(\C),F_1^{\rm pol}(\C))\,\cong\,\left\{\begin{array}{ccc} \C^2 &{\rm for}&p=2 \\
\C &{\rm for}&p=1,3 \\ 0 &{\rm for}&p\not=1,2,3 \end{array}\right.$$
\end{theo}
From the computations, we can deduce more precisely that we can take as generators for $H^p(sl_2(\C),F_0^{\rm pol}(\C))$ for $p=1,2$ (the generator for $p=0$ is the constant function):
$$\theta_1(f(x)\frac{d}{dx}) = f(x)(dx)^0,$$
$$\theta_2(f(x)\frac{d}{dx}) = f'(x)(dx)^0$$
and
\begin{displaymath}
\eta(f(x)\frac{d}{dx},g(x)\frac{d}{dx}) =
\left|\begin{array}{cc} f & g \\ f' & g' \end{array}\right|(x)(dx)^0.
\end{displaymath}
and the generators for $H^p(sl_2(\C),F_1^{\rm pol}(\C))$ for $p=1,2,3$ are 
\begin{displaymath}
\zeta(f(x)\frac{d}{dx}) = f''(dx)^1,
\end{displaymath}
\begin{displaymath}
\omega_1(f(x)\frac{d}{dx},g(x)\frac{d}{dx}) =
\left|\begin{array}{cc} f' & g' \\ f'' & g'' \end{array}\right|(x)(dx)^1,
\end{displaymath}
\begin{displaymath}
\omega_2(f(x)\frac{d}{dx},g(x)\frac{d}{dx}) =
\left|\begin{array}{cc} f & g \\ f'' & g'' \end{array}\right|(x)(dx)^1,
\end{displaymath}
and
\begin{displaymath}
\theta_3(x)(f(x)\frac{d}{dx},g(x)\frac{d}{dx},h(x)\frac{d}{dx}) =
\left|\begin{array}{ccc} f & g & h \\ f' & g' & h' \\ f'' & g'' & h'' \end{array}\right|(x)(dx)^1.
\end{displaymath}
This theorem can easily be computed directly by reducing first the complex to the subcomplex of cochains which are invariant under the Cartan subalgebra of $sl_2(\C)$ (cf theorem 1.5.2 in \cite{Fuk}). 

Let us also recall that $H^3(sl_2(\C),\C)$, $H^3(W_1^{\rm pol}(\C),\C)$ and $H^3(Vect(S^1),\R)$ are all one dimensional and generated by the Godbillon-Vey cocycle (for $sl_2(\C)$, the Godbillon-Vey cocycle is a non-zero multiple of the standard cocycle $\langle [x,y],z\ra$ built from the bracket $[,]$ and the Killing form $\langle,\ra$, see lemma $4$).

Writing down for $sl_2(\C)\subset W_1^{\rm pol}(\C)$ the long exact sequence corresponding to the coefficient sequence from the example of section 3.1, we find
$$0\to H^2(sl_2(\C),F_0^{\rm pol}(\C))\cong\C \stackrel{d_{DR}}{\to} H^2(sl_2(\C),F_1^{\rm pol}(\C))\cong\C^2 \to H^3(sl_2(\C),\C)\cong\C \to 0.$$
It is clear that $d_{DR}$ sends $\eta$ to $\omega_2$, thus $\partial$ must send $\omega_1$ to $\theta(0)$ (with our choices of generators and up to a non-zero multiple). This explains once again the construction of a crossed module for $sl_2(\C)$ representing the generator of $H^3(sl_2(\C),\C)$, cf remark $6$. 

Writing down for $W_1^{\rm pol}(\C)$ the long exact sequence corresponding to the coefficient sequence from the example of section 3.1, we find 
$$0\to H^2(W_1^{\rm pol}(\C),F_1^{\rm pol}(\C))\cong\C\to H^3(W_1^{\rm pol}(\C),\C)\cong\C\to 0.$$
This means that (up to a non-zero multiple) the generator of $H^2(W_1^{\rm pol}(\C),F_1^{\rm pol}(\C))$ is sent to $\theta(0)$ under $\partial$. Thus one can show lemma $1$ by pure cohomological results; it remains to show that $\omega_1$ generates $H^2(W_1^{\rm pol}(\C),F_1^{\rm pol}(\C))$. But this is done in the following lemma:
\begin{lem}
The $2$-cocycle $\omega_1$ can be taken as generator of $H^2(W_1^{\rm pol}(\C),F_1^{\rm pol}(\C))$.
\end{lem}
\pr
Let us denote for the sake of clearness of the proof the generator of $H^2(sl_2(\C),F_1^{\rm pol}(\C))$ by $\omega_1'$ and the cocycle in $C^2(W_1^{\rm pol}(\C),F_1^{\rm pol}(\C))$ having the same formula by $\omega_1$. In the same way, denote by $\theta(0)'$ the (Godbillon-Vey cocycle seen as a) generator of $H^3(sl_2(\C),\C)$, and by $\theta(0)$ the cocycle in $C^3(W_1^{\rm pol}(\C),\C)$ having the same formula. By functoriality with respect to the inclusion $i:sl_2(\C)\subset W_1^{\rm pol}(\C)$, we have a commuting diagram:

\hspace{2.5cm}\xymatrix{
H^2(W_1^{\rm pol}(\C),F_1^{\rm pol}(\C)) \ar[r]^{\partial} \ar[d]^{i^*} & H^3(W_1^{\rm pol}(\C),\C) \ar[d]^{i^*} \\
H^2(sl_2(\C),F_1^{\rm pol}(\C)) \ar[r]^{\partial} & H^3(sl_2(\C),\C)}\vspace{.3cm} 

But it is clear that $i^*(\omega_1)=\omega_1'$, and $i^*(\theta(0))=\theta(0)'$, thus $\partial\omega_1$ must be a non-zero multiple of $\theta(0)$.

The formula $i^*(\omega_1)=\omega_1'$ can also be understood on the level of extensions: the pull-back of the non-trivial extension $\mathfrak{e}$ in $H^2(W_1^{\rm pol}(\C),F_1^{\rm pol}(\C))$ to $sl_2(\C)$ must give $\omega_1'$ (with the associated extension $\mathfrak{e}^{\omega_1'}$) by the above. But there are not many homomorphisms of Lie algebra extensions $\widetilde{\phi}:\mathfrak{e}^{\omega_1'}\to\mathfrak{e}$. By the commutativity and homomorphy constraints, the restriction of $\widetilde{\phi}$, $\phi:F_1^{\rm pol}(\C)\to F_1^{\rm pol}(\C)$ must be a module homomorphism. But $F_1^{\rm pol}(\C)$ is irreducible, thus (as $\phi\not=0$), it must be an isomorphism. But it is easy to see that $F_1^{\rm pol}(\C)$ does not have any non-trivial automorphism ($\phi$ must preserve degrees !). Thus the generator of $H^2(W_1^{\rm pol}(\C),F_1^{\rm pol}(\C))$ must be $\omega_1$.\fin

Observe that we can prove that our construction of a crossed module for $W_1^{\rm pol}(\C)$ representing $\theta(0)$ based on $\omega_1$ is true {\it without using the long exact sequence in continuous cohomology}. Indeed, the crossed module works for $sl_2(\C)$ (ordinary cohomology) and the pullback argument establishes that even for $W_1^{\rm pol}(\C)$, the connecting homomorphism $\partial$ from ordinary cohomology must send $\omega_1$ to $\theta(0)$.

One could try the same discussion for $\mathfrak{g}=Vect(S^1)$, but the long exact sequence induced by the short exact sequence of modules has far too few holes to be accessible. So, it is difficult to follow the cocycles through the sequence, and to see that $\partial c_1=\theta$. But for $Vect(S^1)$, we indicated in \cite{Wag} a method of constructing a crossed module representing $\theta$ based on the result for $W_1$, cf section $5$.

For the purposes of section $5$, we recall the following theorem of Ovsienko- Roger \cite{OvsRog}. Let $\underline{\cal F}_{\lambda}$ be the $Diff(S^1)$-modules of tensor densities on $S^1$. The following maps
$$l:\Phi\mapsto log(\Phi'(x))$$
$$dl:\Phi\mapsto dlog(\Phi'(x))=\frac{\Phi''}{\Phi'}dx$$
$$S:\Phi\mapsto\left[\frac{\Phi'''}{\Phi'}-\frac{3}{2}\left(\frac{\Phi''}{\Phi'}\right)^2\right](dx)^2$$
define 1-cocycles on $Diff(S^1)$ with values in $\underline{\cal F}_0$, $\underline{\cal F}_1$ and $\underline{\cal F}_2$ respectively.

\begin{theo}
The cohomology groups $H^2(Diff(S^1);\underline{\cal F}_{\lambda})$ are non zero only for $\lambda=0,1,2,5,7$, and generated by the following non-trivial cocycles:
\begin{displaymath}
B_0(\Phi,\Psi) = const.(\Phi,\Psi) = \mu(\Phi,\Psi)
\end{displaymath}

\begin{displaymath}
B_1(\Phi,\Psi) = \Psi^*(l\Phi)\,\,dl\Phi
\end{displaymath}

\begin{displaymath}
B_2(\Phi,\Psi) = \Psi^*(l\Phi)\,\,S\Phi
\end{displaymath}

\begin{displaymath}
B_5(\Phi,\Psi) =
\left|\begin{array}{cc}  \Psi^*\,\,S\Phi & S\Psi \\
                        (\Psi^*\,\,S\Phi)' & (S\Psi)'\end{array}\right|
\end{displaymath}

\begin{displaymath}
B_7(\Phi,\Psi) =
2\left|\begin{array}{cc} \Psi^*\,\,S\Phi & S\Psi \\
                        (\Psi^*\,\,S\Phi)''' & (S\Psi)''''\end{array}\right| - 9\left|\begin{array}{cc}(\Psi^*\,\,S\Phi)' & (S\Psi)' \\
                        (\Psi^*\,\,S\Phi)'' & (S\Psi)''\end{array}\right| - \frac{32}{3}(S\Psi + S(\Phi\circ\Psi))B_5(\Phi,\Psi)
\end{displaymath}
\end{theo}
Here $\mu(\Phi,\Psi)$ is the Bott-(Thurston-Virasoro) cocycle which is constant valued.

\end{appendix}

\end{document}